\documentclass[a4paper]{amsart}

\newcommand*{\mykeywords}{Exponential rings, exponential polynomials, exponential ideals}

\title{Note on radical and prime E-ideals}

\author[A. Fornasiero]{Antongiulio Fornasiero}
\address{Universit\`a di Firenze}  
\email{antongiulio.fornasiero@gmail.com} 
\urladdr{https://sites.google.com/site/antongiuliofornasiero/}

\author[G. Terzo]{Giuseppina Terzo}
\address{Universit\`a degli Studi di Napoli "Federico II"}
\email{giuseppina.terzo@unina.it}

\date{13/06/23}

\usepackage{ifpdf}
\ifpdf
\usepackage[pdftex,
pdfborder={0 0 0}, plainpages=false,%
pdfpagelabels, pdfdisplaydoctitle=true,%
pdfstartview={XYZ null null null},%
colorlinks=false
]{hyperref}
\hypersetup{%
  pdfauthor={A. Fornasiero, G. Terzo},
  pdfkeywords={\mykeywords},
  pdfsubject={Exponential ideals in the ring of exponential polynomials:
    maximal and prime}
}
\else
\usepackage[plainpages=false, linktocpage=true, pdfpagelabels, colorlinks=false
]{hyperref}
\fi

\usepackage{amsmath, amsthm, amssymb}
\usepackage{mathtools}
\usepackage{xspace}

\usepackage{comment}

\usepackage[neveradjust]{paralist}
\AtBeginDocument{\setdefaultleftmargin{0pt}{}{}{}{}{}}
\setdefaultenum{\upshape(1)}{\upshape i)}{}{}

\usepackage[T1]{fontenc}

\usepackage[ alphabetic, msc-links]{amsrefs}

\newcommand{\Ffam}{\mathcal F}
\newcommand{\Tbar}{\overline T}
\DeclareMathOperator{\Erad}{E-rad}
\newcommand*{\model}[1]{\mathcal M(#1)}
\newcommand*{\radicale}[1]{\operatorname{{#1}-rad}}
\DeclareMathOperator{\Tradop}{T-rad}
\DeclareMathOperator{\Tzradop}{T_2-rad}
\DeclareMathOperator{\Tunoradop}{T_1-rad}
\newcommand*{\Trad}[1]{\Tradop(#1)}
\newcommand*{\Tzrad}[1]{\Tzradop(#1)}
\newcommand*{\Tunorad}[1]{\Tunoradop(#1)}

\newcommand*{\set}[1]{\{#1\}}

\newcommand{\N}{\mathbb{N}}
\newcommand{\Z}{\mathbb{Z}}

\newcommand{\C}{\mathbb{C}}
\newcommand{\Q}{\mathbb{Q}}

\newcommand{\cv}{\bar c}

\newcommand{\tv}{\bar t}
\newcommand{\x}{\bar x}
\newcommand{\y}{\bar y}

\def\Ind#1#2{#1\setbox0=\hbox{$#1x$}\kern\wd0\hbox to
  0pt{\hss$#1\mid$\hss}\lower.9\ht0\hbox to 0pt{\hss$#1\smile$\hss}\kern\wd0}

\newcommand{\wloG}{w.l.o.g\mbox{.}\xspace}

\newcommand{\ie}{i.e\mbox{.}\xspace}

\newtheorem{lemma}{Lemma}[section]
\newtheorem{thm}[lemma]{Theorem}
\newtheorem{corollary}[lemma]{Corollary}

\newtheorem{proposition}[lemma]{Proposition}
\newtheorem{open problem}[lemma]{Open problem}

\newtheorem*{fact*}{Fact}

\theoremstyle{remark}

\newtheorem{claim}{Claim}
\newtheorem*{claim*}{Claim}

\theoremstyle{definition}
\newtheorem{definition}[lemma]{Definition}

\newtheorem{remark}[lemma]{Remark}
\newtheorem{final remark}[lemma]{Final remark}
\newtheorem{example}[lemma]{Example}
\newtheorem{examples}[lemma]{Examples}

\begin{document}

\begin{abstract}
We show that $\C[\x]^{E}$ is not Noetherian
even respect to prime E-ideals. Moreover we give a characterization of exponential radical ideals 
\end{abstract}

\maketitle


\makeatletter
\renewcommand\@makefnmark%
   {\normalfont(\@textsuperscript{\normalfont\@thefnmark})}
\renewcommand\@makefntext[1]%
   {\noindent\makebox[1.8em][r]{\@makefnmark\ }#1}
\makeatother

\section{Introduction}

The notion of exponential ideal (E-ideal) was introduced in the papers
\cites{van, macintyre}; it was related to the study of exponential function
after the problem posed by Tarski on the decidability of the reals with
exponentiation. 
An exponential ring (E-ring) is a pair $(R, E)$ where R is a
commutative ring with 1 and E is a homomorphism
$E: (R, +) \rightarrow (R^{*}, \cdot)$.
We will always assume that $R$ is a $\Q$-algebra.
The classical examples are the reals and the complex
numbers with the usual exponentiation. 

Starting from an E-ring $R$,
we can construct the E-polynomial ring in the variables $\x$ over $R$ 
by induction (see
\cites{van, macintyre}) and we denote it by $R[\x]^E$. 
An E-ideal $I$ of an E-ring $R$ is an ideal with the property that if $\alpha \in I$ then
$E(\alpha) - 1 \in I$; they coincide with the possible kernels of homomorphisms which
preserve the exponential map~$E$. 

Some further contributions to the study of E-ideals
were given in \cites{manders, point, terzo}.  
In the recent paper \cite{ideali}, together with P. D'Aquino, 
we studied E-ideals of E-rings,  
we gave two notions of maximality for E-ideals, 
and related them to primeness. 
We proved that the three notions are independent, unlike in the classical
case. Moreover, we showed that, for any exponential field $K$,
not all maximal E-ideal of $K[\x]$ correspond to 
points of $K^n$. 

In this paper we further investigate E-ideals of $R[\x]^{E}$.
It was known that the exponential Zariski topology on $\C^{n}$ is not Noehterian
(see \cite{macintyre}).
We show that $\C[\x]^{E}$ is not Noetherian
even respect to prime E-ideals. 
Moreover, we give a reasonable notion of E-radical
E-ideals, characterize them, and prove some of their properties, using
a technique introduced in \cite{ideali} that permits to extend prime ideals to
prime E-ideals.

\section{E-polynomial ring and basic results}

First to introduce the construction of the E-polynomial ring is useful to recall the notion of partial E-rings and partial E-ideals and some important result related to it proved in \cite{ideali}.

\subsection{Partial E-rings and E-ideals}

\begin{definition}
A partial E-ring is a triple $D=(D, V, E)$ where
\begin{enumerate}
\item $D$ is a   $\Q$-algebra; 
\item $V $ is a  $\Q$-vector subspace of
  $D$ containing~$\Q$;
\item $E : (V, +) \to (D^{*}, \cdot)$ is a group homomorphism.
\end{enumerate}
\end{definition}

\begin{definition}
A partial E-ideal of $D$ is an ideal $I$ of the ring $D$ such that, for every
$v \in I \cap V$, $E(v) - 1 \in I$. If $D=V$ we say that it is an E-ideal.
\end{definition}

\begin{remark}
\label{trivialEideal}
If $I$ is an ideal of $D$ with $I\cap V=( 0)$ then $I$  is an E-ideal of $D$.
\end{remark}

\begin{definition}
An E-ideal I of D is prime if it is prime as ideal, i.e. iff $D/I$ is a domain. An E-ideal I of D is an E-maximal ideal if it is maximal among the E-ideals. It is strongly maximal if it is maximal
as ideal.
\end{definition}

\subsection{Construction of E-polynomial ring}

The construction of the E-polynomial ring in many variables is well known, see
\cites{van, macintyre};
for the reader's convenience, we briefly recall it here. 
Starting from an E-ring $R$, we construct  the E-polynomial ring in the
variables $\x = (x_1, \ldots, x_n)$, denoted by $R[\x]^E$, as a union of a chain of partial E-rings equipped with partial E-morphisms. 
The following three chains are constructed by recursion: 
 $(R_{\scriptsize k}, +, \cdot)_{\scriptsize {k\geq-1}}$, $(B_k, +)_{k\geq0}$, and  $(E_k)_{k\geq-1}$ will be rings, abelian groups and partial E-morphisms, respectively. 

Let $R_{-1} = R$, and $R_0 = R[\overline x]$ as partial exponential ring (where
$\exp$ is defined only on~$R$).
Let $A_0 = (\overline x) $ be the ideal
of $R[\overline x]$ generated by $\overline x$. 
For $k \geq 1$, let $A_{k}$
be the $R_{k-1}$-submodule of $A_{k}$  
generated by $t^a$ with  $0 \neq a \in A_k$.
So, we define  $R_k = R_{k-1} \oplus A_k,$ $R_{k+1} = R_{k}[t^{A_k}]$ and $E_{k} : (R_k, +) \rightarrow (R_{k+1}^{*}, \cdot)$ such that $E_k(x) = E_{k-1} (r) \cdot t^b$, for $x = r + b$, $r \in R_{k - 1}$ and $b \in B_k$.
Then the E-polynomial ring is the limit of this object, i.e $R[\overline x]^E= \cup R_k.$
Sometimes it is convenient to represents $R[\x]^{E}$ as the group ring  $R[\x][t^{\bigoplus_{i\geq 0} A_i}]$. 

In \cite{ideali} we generalized this construction to any partial E-ring R, \ie we gave a free completion of a partial E-ring R that they denoted by $R^E$. Moreover we gave sufficient conditions on a subring S of $R^E$ such that the free completion of $S,$ $S^E$ is isomorphic to $R^E.$ We recall the result just for E-polynomial rings.

\begin{lemma}[\cite{ideali}*{Lemma~2.10}] \label{completion}
Let $R$ be an E-ring and $S$ be a partial subring of $R[\x]^E,$ and assume that
$S = R[\x][e^A]$ for some $\Q$-linear subspace $A$ of $R[\x]$ which has trivial
intersection with~$R$.
Then, $S^E = R[\x]^E$.
\end{lemma}

\begin{lemma} \label{idealiprimi}
Let $S \subseteq R[\x]^{E}$ be as in Lemma~\ref{completion}.
Let $I$ be an ideal of  $S$.
If $I$ is a partial prime E-ideal of $R_n$, then $I^E$ (the $E$-ideal of
$R[\x]^{E}$ generated by $I$) is a prime E-ideal of $R[\x]^E$.
\end{lemma}

\section{Noetherianity}

As Macintyre point out in  \cite{macintyre},  
neither $\Bbb C[x]^E$ nor $\Bbb R[x]^E$ are Noetherian for E-ideals, by considering the E-ideal $I = (E(\frac{x}{2^n}) - 1)_{n \in \Bbb N},$ since it is not finitely generated (for details see also \cite{terzo}).

\smallskip

There is a notion of noetherianity also for topological spaces:
\begin{definition}
A  topological space $X$ is Noetherian if it satisfies the descending chain condition for closed subsets, \ie any strictly descending sequence of elements of $X$ is stationary. \end{definition}

$\C[\x]^{E}$ is very far from being Noetherian for E-ideals, since
$\C$ with the exponential Zariski topology is not a Noetherian space, because of
the chain:
\begin{align*}
Z(E(x) -1) & \supset  Z(E(\frac{x}{2}) -1) \supset  \ldots \supset  Z(E(\frac{x}{n!}) -1)\supset  \ldots\\
\{2\pi i \Z\} &\supset \{4\pi i \Z\} \supset \ldots \supset \{2 n! \pi i \Z\} \supset \ldots 
\end{align*}

descendes indefinitely.
 
\smallskip

We show that $\C [\x]^E$ is not Noetherian even respect to prime E-ideals.

\begin{thm}
The ring $\C [\x]^E$ doesn't satisfy ACC condition for prime E-ideals.

\end{thm}

\proof
We consider the subring $S := \C[\x, e^{\C \x}] =  \C[\x] [e^{\C\x}]  $
of the ring $\C [\x]^E$. 
The idea is to construct an ascending chain of prime
E-ideals of $S$ and to extend it to be an ascending chain of prime E-ideals of  
$\C [\x]^E$.
Let $P_i < S$ be prime ideals of $S$ and we define $Q_i = P_i  \C [\x]^E.$

For simplicity, we consider the case when we have only one variable 
$\x = x$.
We construct prime ideals in the following way: let be $B$ a transcendence basis of $\C$; $B = (b_j : j < 2^{\aleph_0}).$ For all $i \in \Bbb N,$ we define $$p_i := e^{b_{3i}x} +e^{b_{3i + 1}x} + e^{b_{3i +2}x}  .$$
Let be $A_n = (p_0, \ldots, p_n)$ as an ideal in $\C[e^{\C x}].$ 
We need the following result:

\begin{lemma}
The ideal $A_n$ is prime for all $n \in \Bbb N.$
\end{lemma}

\begin{proof}
We introduce  new variables denoting elements of the form $e^{b_i x}$ for any
$b_i \in B$, \ie  $z_i = e^{b_i x}$, so we can denote 
$A_n = (z_0 + z_1 + z_2, z_3 + z_4 + z_5, \ldots, z_{3n} + z_{3n+1} + z_{3n + 2})$ ad
an ideal of $\C [\overline z^{\Q}].$
We prove by induction on $n$ that $A_n$ is prime. For $n=1$ we have that $A_1 =
(z_0+z_1+z_2)$ as an ideal of $\Bbb C [z_0^{\Bbb Q}, z_1^{\Bbb Q}, z_2^{\Q}]$. 
Assume that $p(z_0, z_1, z_2)\cdot q(z_0, z_1, z_2) \in A_1$ then $p(z_0, z_1, z_2)\cdot
q(z_0, z_1, z_2) = r(z_0, z_1, z_2) (z_0+z_1+z_2)$; let $k$ be the common
denominator of any exponents in $p, q, r$, so we can consider $p, q, r \in \Bbb C
[z_0^{\pm \frac{1}{k}}, z_1^{\pm \frac{1}{k}}, z_2^{\pm \frac{1}{k}}]$. 
By replacing $z_{i}$ with $t_{i}^{k}$, we have that
$p,q,r \in \C[t_{0}^{\pm 1}, t_{1}^{\pm 1}, t_{2}^{\pm 1}]$.
Notice that $z_{0} + z_{1} + z_{2}$ becomes $s(\tv) := t_{1}^{k} + t_{2}^{ k} +
t_{3}^{k}$, and that $s$ is irreducible in $\C[\tv]$ and hence the ideal
$A_{1}'$ generated by $s$ inside $\C[\tv]$ is prime.
By localization, the ideal $A_{1}''$ generated by $s$ inside $\C[\tv^{\Z}]$ is
also prime; therefore, since $p q \in A_{1}''$, we have either $p \in A_{1}'' \subseteq
A_{1}$ or $q \in A_{1}'' \subseteq A_{1}$, proving that 
$A_1$ is prime.

For $n = 2$, $A_2 = (z_0+z_1+z_2, z_3+z_4+z_5)$ this is prime since 
\[
\frac{\C [z_0^{\Q}, \ldots, z_5^{\Q}]}{A_2} \cong \frac{\C [z_0^{\Q}, z_1^{\Q}, z_2^{\Q}]}{A_1} \otimes
\frac{\C [z_3^{\Q}, z_4^{\Q}, z_5^{\Q}]}{(p_1)}
\]

So $\frac{\Bbb C [z_0^{\Bbb Q}, \ldots, z_5^{\Bbb Q}]}{A_2}$ is a domain, since the
tensor product of $\C$-algebrae which are domains is also a domain, see
\cite{bourbaki}*{Chapter V, \S17}. 
In a similar way we can prove that $A_n$ is prime for any $n \in \Bbb N.$
\end{proof}

\smallskip

If $A_n$ is prime then $P_n = A_n S$ is partial prime E-ideal in S and so $Q_n$, by Lemma \ref{completion} and Lemma \ref{idealiprimi}, will be prime E-ideals in  $\Bbb C [\x]^E.$ So we have an ascending chain of prime E-ideals.
\qed

\bigskip

Now we  give conditions to say when an E-ideal of a particular form
is prime.
Let $\x, \y$ be tuples of variables of the same lenght.
Given $p(\x,\y) \in \C[\x,\y]$, we denote by $\tilde p(\x) := p(\x, E(\x)) \in
\C[\x]^{E}$.
Let $I \subseteq \C[\x,\y]$ be an ideal, and let $\tilde I := \set{\tilde p: p \in I}$.
We denote by $S = \Bbb C[\overline x, e^{\overline Q \overline x}]$,
by $J$ the E-ideal of $\C[\x]^{E}$ generated by $\tilde I$,  and by $H = J \cap S$. 

\begin{proposition} $J$ is prime in $\Bbb C[\overline x]^E$ iff $H$ is prime in $S.$
\end{proposition}

\begin{proof}
One direction is trivial. For the other one we assume that $H$ is a prime in $S$ then by Lemma \ref{completion} and Lemma \ref{idealiprimi} we have that $J$ is a prime E-ideal of $\Bbb C[\overline x]^E.$
\end{proof}

\begin{proposition}
If the following hold:
\begin{enumerate}
\item I is a prime ideal;
\item I doesn't contain non zeros elements of the form $a + \overline q \cdot \overline x$ with $\overline q \in \Bbb Q^n$ and $a \in \Bbb C;$
\item I doesn't contain any
element of the form $\overline y^{\overline q} - a$ where  $\overline q \in \Bbb Q^n$ and $a \in \Bbb C;$
\item For all $n \in \N$  $I_n$ the ideal in $\C[\x, \y^{\frac{1}{n}}]$ is prime.
\end{enumerate}
Then J is a prime E-ideal. 
\end{proposition}

\begin{proof}
We denote by $K$ the ideal generated by $\tilde I$ in $S.$ If we prove that $K$ is a partial E-ideal of S and it is prime we conclude the proof, because by Corollary 3.13 of \cite{ideali} we obtain that $J$ is prime and $K = H.$
First we prove that $K$ is a partial  E-ideal of $S.$ We note that the domain of the exponential map on $S$ is $\Q \cdot \x + \C,$ so we have to prove that $K \cap  (\Q \cdot \x + \C )= (0).$ By replacing $e^{\Q\x}$ with $\y^{\Q}$ we have $\C[\x, \y] \subset \C[\x, \y^{\Q}]\cong \C[\x, e{\Q\x}] = S.$ We claim that $$K \cap \C[\x, \y]  =I.$$ If $a \in K \cap \C[\x, \y] $ then $a = \frac{p}{\y^{m}}$ where $p \in I.$ Then $a\cdot\y^m \in I,$ but by (2)+(3) we have that $a \in I$ and this implies that $K$ is a partial E-ideal.\
Now we have to prove that $K$ is prime. We are assuming that the ideal $I$ of $\C[\x, \y]$ is prime so by commutative algebra (see \cite{Matsumura}) we know that the ideal $I'$ in the localization $ \C[\x, \y^{\pm 1}]$ is prime or trivial, but it is prime because we are assuming (3). If we consider the ideal $I'' = K$ in the integral extension of $\C[\x, \y^{\pm 1}]$, i.e. in the extension  $ \C[\x, \y^{\Q_{\geq 0}}].$ We first observe that $I'' = \cup_{n} I_n,$ from (4) any $I_n$ is prime and this implies that $I''$ is prime, this conclude the proof.

\end{proof}

\section{Exponential radical ideals}

We give the notion of E-radical ideal for any E-ring $R$ as follows.

\begin{definition} Let $J$ an E-ideal of an E-ring $R$. 
We define the E-radical ideal of $J$ as $\Erad(J) :=  \cap_{P \supseteq J} P$, where $P$
varies among prime E-ideals. 
\end{definition}

Let  $J$ an E-ideal of the E-polynomial ring $K[\x]^E$.
Let be $F$ an E-ring containing $K$. 
 
We can define  $\mathcal{I} (V(J))$ as follows:\\
$V_F(J) := \{ \overline a \in F^n : f(\overline a) = 0 \mbox{ for all } f(\x) \in J
\}$,\\
$V(J) := \bigcup V_F(J)$ as $F$ varies among all E-fields containing $K$, and\\
$ \mathcal{I} (V(J)) = \{ p(\x) \in K[\x]^E : p(\overline a) = 0 \mbox{ for all } \overline
a \in V(J) \}$. 

\begin{remark}
Let $(R,E)$ be an E-domain and $K$ be its fraction field.
Then, there exists at least one way to extend the exponential function to all of $K$.
\end{remark}
\begin{proof}
Let $A\subset K$ be a complement of $R$ as $\Q$-linear spaces.
For every $r \in R$ and $a \in A$, define $E'(a + r) := E(r)$.
Then, $E'$ is an exponential function on $K$ extending $E$.
\end{proof}

\begin{corollary}
$\mathcal{I} (V(J)) = 
\{ p(\x) \in K[\x]^E : p(\overline a) = 0 \mbox{ for all } \overline a \in V_{F}(J)$ 
as $F$ varies among all E-domains containing $K\}$. 
\end{corollary}

\begin{lemma}
Let be $J$ an E-ideal of the  E-polynomial ring $K[\x]^E.$ Then $\Erad(J) = \mathcal{I} (V(J)) $.
 \end{lemma}

\begin{remark}
The E-ideal $I = (xy)^E$ is not a E-radical ideal, unlike in the classical
case. Indeed $I$ is not the intersection of prime E-ideal, since $(xy)^E \neq (x)^E
\cap (y)^E,$ because $x(E(y) - 1) \in (x)^E \cap (y)^E \setminus (xy)^E$. 
\end{remark}

\begin{remark}
If $J$ is not contained in any prime E-ideals then $\Erad(J) = K[\x]^E.$ 
We use the technique introduced in  \cite{ideali} to construct such example. Let
be $I = (xy, E(x) + 1, E(y) + 1)$ an ideal of $S = K[x, y, e^{\Bbb Q x}, e^{\Bbb
  Q y}]$ where $S$ is a subring of $K[\x]^E,$ in particular it is a partial
$E$-ideal of $S$. 
$I^E$ is then an E-ideal of $K[\x]^E$ which is not contained in any
prime E-ideal; indeed, if $I \subseteq P$ where $P$ is E-prime, then $xy \in I \subseteq P.$ 
Since $P$ is  prime,  \wloG   $x \in P$, and so $E(x) - 1 \in P$.
But $E(x) + 1 \in I \subseteq P$, and thus $2 \in P$,  contradiction. 
\end{remark}

\section{Characterization of radical E-ideals}
We consider an  E-ring $R$ and let be $J$ an E-ideal of~$R$.  

We study prime E-ideals and  E-radical ideals, i.e. E-ideals which are equal to
their E-radical.
We characterize $\Erad(J)$ using the following theory.\\

Let be $\mathcal{L} = \{+, -, \cdot, e^{x}, 0, 1\} \cup \{ V \}$ where $V$ is a unary predicate.\\

We recall the following definition:


\begin{definition}
Given $p_1, \ldots, p_k, q$ $\mathcal{L} $-terms,  we define the associated strict Horn clause the formulas of the type:
$$V(p_1) \wedge \ldots \wedge V(p_k) \rightarrow V(q).$$ 
\end{definition}

We denote by $H$ the set of all $\mathcal{L} $-strict Horn clause as above.

\begin{examples}
$V(x^2)  \rightarrow V(x)$ is in $H$.
Also $V(x) \wedge V(y)  \rightarrow V(x + y) ,$ $V(x)  \rightarrow V(e^x - 1)$ and $V(0)$ are in~$H$.
\end{examples}

In order to lighten the notation, we write $p_1 \wedge \ldots \wedge p_k \rightarrow q$ in place of
 $V(p_1) \wedge \ldots \wedge V(p_k) \rightarrow V(q)$.

\smallskip

We consider the following theories:\\
\[\begin{aligned}
T_{1} &:= \{ \mbox{ Horn clause } \alpha : \mbox{ for any E-ring R and E-ideal J, } (R, J) \models \alpha \}\\
T_2 &:=  \{ \mbox{ Horn clause } \alpha : \mbox{ for any E-ring R and prime E-ideal
  J, } (R, J) \models \alpha \}.
\end{aligned}
\]

Clearly, $T_{1}$ is generated by the following Horn clauses:
\[ 0, \quad x \wedge y  \rightarrow x-y, \quad x  \rightarrow xy, \quad x \rightarrow e^{x}-1.\]

Our aim is to give an explicit description of $T_{2}$ and relate it to the E-radical.
In $T_2$ we have, besides the clauses in $T_{1}$, also others; for instance,
the following clauses are in $T_{1}$:
$x^n  \rightarrow x$,  $(xy \wedge e^{x} + 1)  \rightarrow y$.

Let $T$ be  a set of Horn clauses.
Let $\model{T}$ be  following family   of subsets of $R$:
\[
\model{T}  := \set{J \subseteq R: (R, J) \models T} 
\]
\begin{remark}\label{unione e intersezione} 
$\model{T}$ is closed under arbitrary intersections and 
under increasing unions.
\end{remark}
Thus, we can consider the ``radical'' operator associated to $T$.
We will let $X$ vary among subsets of~$R$.

We define
\[
\Trad X \coloneqq \bigcap \model{T}.
\]
We have that $\Trad X$ is the smallest subset of $R$ containing $X$ and such
that $(R, \Trad X ) \models T$.
In particular, $(R, X) \models T$ iff $X = \Trad X$.

\begin{example}
$\Tunorad X = (X)^{E}$, the E-ideal generated by~$X$.
\end{example}

\begin{example}
Let $T_{0} := \set{0, x \wedge y \rightarrow x - y, x  \rightarrow xy}$.
Then, $\model{T_{0}}$ is the family of pairs $(R,J)$ with $J$ ideal of $R$.
We have that $\radicale{T_{0}}(X)$ is the ideal generated by~$X$.

Let $T_{3} := \set{0, x \wedge y \rightarrow x - y, x  \rightarrow xy, x^{2} \rightarrow x}$.
Then, $\model{T_{3}}$ is the family of pairs $(R,J)$ with $J$ radical ideal of $R$.
We have that $\radicale{T_{3}}(X)$ is the radical of the ideal generated by~$X$.
\end{example}

We can build $\Trad J$ in a  ``constructive'' way: \ie, we have a description
of all elements of $\Trad J$.

\begin{definition}
Given a family $\Ffam$ of 
$\mathcal{L}$-structures, we denote its theory by
\[
Th(\Ffam) := \set{\alpha \in H: \forall M \in \Ffam\ M \models \alpha }.
\]
The deductive closure of $T$ is $\Tbar := Th(\model T)$.
We say that  $T$ clauses is deductively closed if $T = Th(\Ffam)$ for some
family $\Ffam$, or equivalently if $T = \Tbar$.
An axiomatization of $T$ is a set of Horn clauses $S$ such that
$\overline S = \Tbar$.
\end{definition}

\begin{remark}
\[
\Trad X  = \set{q(\cv): 
p_{1}(\x) \wedge \dots \wedge p_{k}(\x) \rightarrow q(\x) \in \Tbar,
 \cv \in R^{< \omega}, p_i(\cv) \in X, i = 1, \dotsc, k}. 
\]
\end{remark}
\begin{proof}
Let us denote 
$Y := \set{q(\cv): 
p_{1}(\x) \wedge \dots \wedge p_{k}(\x) \rightarrow q(\x) \in \Tbar,
 \cv \in R^{< \omega}, p_i(\cv) \in X, i = 1, \dotsc, k}$.

It is clear that $X \subseteq Y \subseteq \Trad X$.
We want to show that $\Trad X \subseteq Y$.
It suffices to show that $Y \in \model T$.
Let $p_{1}(\x) \wedge \dots \wedge p_{k}(\x) \rightarrow q(\x) \in T$ and $\cv \in R^{< \omega}$ such
that
$p_{i}(\cv) \in Y$, $i = 1, \dotsc, k$. 
It suffices to show the following:
\begin{claim}
$q(\cv) \in Y$.
\end{claim}
For simplicity of notation, we assume that  $k = 1$ and $p := p_{1}$.
Since $p(\cv) \in Y$, by definition of~$Y$, there exist $\beta := r_{1}(\x, \x') \wedge
\dots \wedge r_{\ell}(\x, \x') \rightarrow p(\x) \in T$ and $\cv' \in R^{< \omega}$
such that $r_{j}(\cv,\cv') \in X$, $i = 1, \dotsc, \ell$.
Notice that every $M \in \model T$ satisfies
\[
\beta := \bigwedge_{\ell} r_{\ell}(\x, \x') \rightarrow q(\x).
\]
Therefore, $\beta \in  \bar T$, and therefore, by definition,
$q(\cv) \in Y$.
\end{proof}

\begin{corollary}
For every $b \in R$,
\[
\Trad{Xb} = \set{q(\cv): 
p_{1}(\x) \wedge \dots \wedge p_{k}(\x) \rightarrow q(\x) \in \Tbar,
 \cv \in R^{< \omega}, p_i(\cv) \in X \vee p_{i}(\cv) = b, i = 1, \dotsc, k}. 
\]
\end{corollary}

The following theorem gives more explicit description of 
$\Tzrad{X}$.

\begin{thm}\label{thm:Erad-up}
For every $n \in N$, we define the operator $\sqrt[n]{\vphantom{x} }$ on subsets of $R$
inductively in the following way:\\
$\sqrt[0] X = (X)^E;$ \\
$\sqrt[1] X =  \sqrt[0] { \{ a \in R:  
\exists b_1, b_2 \in R : b_1 \cdot b_2  \in
  \sqrt[0] X  \wedge a  \in \sqrt[0]{X b_1} \cap \sqrt[0]{X b_2} \} };$ \\
$\sqrt[n+1] X = \sqrt[n] { \{ a :  
\exists
b_1, b_2  : b_1 \cdot b_2 \in  \sqrt[n] X  \wedge a  \in \sqrt[n] {X b_1} \cap \sqrt[n]{X b_2} \} }.$\\
Define $\sqrt[E] X := \bigcup_{n \in \N} \sqrt[n] X$

Then, $\Erad X = \Tzrad X = \sqrt[E] X$.
\end{thm}

\proof
We first need some results.
\begin{lemma}\label{lem:Erad-prod}
Let $J \subseteq R$.
Let $b_{1} \cdot b_{2} \in J$.
Assume that $a \in \Tzrad{J b_{1}} \cap \Tzrad{J b_{2}}$.
Then, $a \in \Tzrad J$.
\end{lemma}
\begin{proof}
If $J$ were a prime E-ideal, the result would be clear.
In general, there exist Horn clauses
$p_{i,1}(\x) \wedge \dots \wedge p_{i,k}(\x) \wedge z \rightarrow q_{i}(z,\x) \in I(T)$, $i = 1,2$,
and
$\cv \in R^{< \omega}$, such that
\[
a = q_{i}(b_{i},\cv), 
\quad
p_{i,j}(\cv) \in J,
\quad i = 1,2, 
\quad j = 1, \dotsc, k.
\]
Moreover, the Horn clause
\[
\alpha(w, z_{1}, z_{2}, \x)  := z_{1} \cdot z_{2} \wedge w - q_{1}(z_{1}, \x) \wedge w - q_{2}(z_{2},\x)
\wedge \bigwedge_{i,k} p_{i,k}(\x) 
\rightarrow w
\]
is in $T_{0}$ (since it is satisfied by any prime E-ideal).
Thus, $(R, \Tzrad J) \models \alpha$ and the conclusion follows by considering
$\alpha(a, b_{1}, b_{2}, \cv)$.
\end{proof}

It is therefore that $\Erad X \supseteq \Tzrad X \supseteq \sqrt[E] X$.
Thus, it suffices to show that $\Erad X \subseteq \sqrt[E] X$, or, equivalently, that
for every $a \in R \setminus \sqrt[E]X$, we have $a \notin \Erad X$.

 



We need a further result.

\begin{lemma}\label{lem:Trad-max}
Let  $a \in R$ and $P$ be an E-ideal of $R$ maximal among the E-ideals $J $of $R$
not  containing~$a$ and such that $\sqrt[E] J = J$.
Then, $P$ is prime.
\end{lemma}

\begin{proof}
Suppose $b_1 \cdot b_2 \in P$ and by contradiction we assume that $b_1, b_2 \not \in P$.
We define $Q_1 = \sqrt[E]{P  b_1}$ and $Q_2 = \sqrt[E]{P  b_2}$. 
We have that  $Q_{1},Q_{2} \in \model{T_{0}}$.
Notice that $Q_1, Q_2 \supset P$, and therefore the maximality of $P$ implies that $a
\in Q_{1}$ and $a \in Q_{2}$.
Moreover, $b_1 \cdot b_2 \in P$; since $\sqrt[E] P = P$, we have
  that $a \in  \Tzrad {P } = P$,  contradiction.   
\end{proof}

We can now conclude the proof of Thm.~\ref{thm:Erad-up}.
Let  $a \notin \sqrt[E]X$.
By Lemma~\ref{lem:Trad-max}, there exists a prime E-ideal $P$
containing $X$ such that $a \notin P$.
Thus, $a \notin \Erad X$, and we are done.
\qed

\begin{corollary}\label{cor:Erad}
$J$ is a $E$-radical iff, for every $a, b_{1}, b_{2} \in R$
\[
b_{1} \cdot b_{2} \in J \wedge a \in \sqrt[E]{J b_{1}} \cap \sqrt[E]{J b_{1}} 
\rightarrow a \in J.
\]
\end{corollary}
From the above Corollary we can extract a recursive axiomatization of $T_{0}$, 
using
Thm.~\ref{thm:Erad-up} to characterize $\sqrt[E]{J b_{i}}$.

\bigskip

We can interpret the discussion above in terms of quasi-varieties.
\begin{definition}
An E-ring is E-reduced if $\Erad(0) = (0)$.
\end{definition}
By Thm.~\ref{thm:Erad-up}, the class E-red
of E-reduced E-rings is a quasi-variety:
\ie, it can be axiomatized via Horn formulae in the language
\[
(=, + , - ,\cdot, e, 0, 1)
\]
(we replace the condition $t \in (0)$ with the condition $t = 0$).
By \cite{Burris-Sankappanavar}*{Thm.~2.25}, E-red is closed under isomorphisms, taking substructures, reduced product.\\
It is not difficult to see directly that E-red is closed under isomorphisms,  substructures, direct
products,  and ultraproducts: thus, by \cite{Burris-Sankappanavar}*{Thm.~2.25},
E-red is a quasi-variety.
From this it is easy to deduce that $\Erad = \Tzradop$.
With the proof we gave of Thm.~\ref{thm:Erad-up} we obtained the extra
information that $\Erad = \sqrt[E]{\vphantom{x}}$, and Corollary~\ref{cor:Erad},
with the recursive axiomatization of~$T_{0}$.

\bibliographystyle{plain}
\bibliography{exponential}		

\begin{thebibliography}{99}

\bibitem{Burris-Sankappanavar}
S.Burris, H. P. Sankappanavar:
\emph{A Course in Universal Algebra}, The Millennium Edition.
S. Burris and H.P. Sankappanavar \url{https://www.math.uwaterloo.ca/~snburris/htdocs/ualg.html}

\bibitem{bourbaki} N. Bourbaki:  \emph{Elements of Mathematics- Commutative Algebra: Chapters 1-7}, Springer.

\bibitem{ideali} P. D'Aquino, A. Fornasiero and G. Terzo: \emph{E-ideals in exponential polynomial ring}, Submitted.
 
\bibitem{pag} P. D'Aquino, A. Macintyre and G. Terzo: \emph{Schanuel Nullstellensatz for Zilber
fields}, Fundamenta Mathematicae 207, 123-143 (2010).

\bibitem{shapiro} P. D'Aquino, A. Macintyre and G. Terzo: \emph{From Schanuel Conjecture to Shapiro Conjecture}, Commentarii Mathematicae Elvetici 89 (3), 597-616 (2014).



\bibitem{van} L. van den Dries: \emph{Exponential rings, exponential polynomials and exponential functions}, Pacific Journal of Mathematics, 113, (1), 51-66 (1984).

\bibitem{Everest_Porten} G. R. Everest and   A. J. van der Poorten:
\emph{Factorisation in the ring of exponential polynomials}, of
Proceedings of the American Mathematical Society, 125, (5),
(1997), 1293-1298.


\bibitem{HR} Henson,  Rubel and Singer:  \emph{Algebraic Properties of the Ring of General Exponential Polynomials} Complex Variables, 13,  1-20 (1989).

\bibitem{macintyre} A. Macintyre:  \emph{Exponential Algebra}  In: Ursini, A., et al., eds. Logic and Algebra. Proceedings of the international conference dedicated to the memory of Roberto
Magari. Lect. Notes Pure Appl. Math 180. M. Dekker, pp. 191–210.

\bibitem{manders} K. Manders:  \emph{On algebraic geometry over rings with exponentiation} Z. Math. Logik Grundlag. Math. 33 (1987), no. 4, 289-292.


\bibitem{marker}
D. Marker: \emph{A remark on Zilber's pseudoexponentiation}, The Journal of Symbolic
Logic, 71, (3),  791-798 (2006).


\bibitem{Matsumura}
H. Matsumura: \emph{Commutative Algebra}, W. A. Benjamin, Inc., New York, 1970.

\bibitem{point}
F. Point and N. Regnault  \emph{Exponential ideals and a Nullstellensatz}, arXiv:2004.10444, 22 April 2020


\bibitem{Ritt} J.F. Ritt: \emph{A factorization theorem of functions $\sum_{i=1}^na_ie^{\alpha_iz}$}, Transactions of
American Mathematical Society 29, (1927), 584-596.


\bibitem{terzo} G. Terzo: \emph{Some consequences of Schanuel's conjecture in exponential rings.} Communications in  Algebra 36, 3, (2008), 1171-1189.

\end{thebibliography}

\end{document}